\documentclass{article}

\usepackage[T2A]{fontenc}
\usepackage[cp1251]{inputenc}
\usepackage[english,russian]{babel}
\usepackage[tbtags]{amsmath}
\usepackage{amsfonts,amssymb}

\usepackage{cite}

\usepackage{mathrsfs}
\usepackage[short]{mz2ewin}
\usepackage{graphicx}

\theoremstyle{plain}

\newtheorem{lemma}{Лемма}
\newtheorem{theorem}{Теорема}

\theoremstyle{definition}

\newtheorem{remark}{Замечание}

\begin{document}

\udk{517.968}


\author{С.\,А.~Бутерин}
\address{Саратовский национальный исследовательский университет имени Н.Г. Чернышевского;\\
Московский центр фундаментальной и прикладной математики;\\
Московский государственный университет имени М.В. Ломоносова}
\email{buterinsa@sgu.ru}

\title{Об управляемой системе на бесконечном  временн\'ом дереве}

\maketitle
\begin{fulltext}




\footnotetext{Работа выполнена при финансовой поддержке РНФ (проект \No\ 24-71-10003).}

Н.Н.~Красовским в~\cite{Kras-68} была поставлена и исследована задача об успокоении управляемой системы с
последействием для уравнения первого порядка запаздывающего типа с постоянными вещественными коэффициентами
на интервале. Позднее А.Л.~Скубачевский~\cite{Skub-94} рассмотрел более трудный случай уравнения нейтрального
типа
\begin{equation}\label{1}
y'(t)+ay'(t-\tau)+by(t)+cy(t-\tau)=u(t), \quad t>0, \qquad a,b,c\in{\mathbb R}, \; \tau>0.
\end{equation}
Предыстория системы определяется условием $y(t)=\varphi(t)\in W_2^1[-\tau,0],$ $t\in[-\tau,0].$ Для
фиксированного $T>2\tau$ требуется подобрать управление $u(t)\in L_2(0,T),$ приводящее систему в состояние
равновесия $y(t)=0$ при $t\in[T-\tau,T].$ Были установлены существование и единственность оптимального
управления $u(t),$ т.е. минимизирующего $\|u\|_{L_2(0,T)}.$ При этом оптимальная траектория $y(t)$
удовлетворяет некоторой краевой задаче для уравнения второго порядка с разнонаправленными сдвигами аргумента,
решение которой понимается в обобщенном смысле, коль скоро $a\ne0$ (также см. \cite{Skub-97}).  Дальнейшие
обобщения данной задачи на интервале можно найти в \cite{AdSkub-24} и приведенных там ссылках.

В \cite{But23-arXiv, But24} при помощи идеи глобального запаздывания \cite{But23-RM} данная задача перенесена на так называемые квантовые
графы и был охвачен случай произвольного конечного дерева. Это привело к концепции временн\'ого графа, когда  переменная на ребрах
отождествляется со временем, а каждая внутренняя вершина понимается как точка разветвления процесса на несколько параллельных процессов по
числу выходящих из нее ребер.

Систему управления на дереве можно интерпретировать следующим образом. В каждый момент времени, соответствующий какой-либо внутренней
вершине, появляется не\-сколько вариантов дальнейшего течения процесса. Требуется подобрать управление, приводящее систему в состояние
равновесия, независимо от того, какой именно набор сценариев в итоге будет реализован. Было установлено, что оптимальная траектория
дополнительно удовлетворяет условиям типа Кирхгофа во внутренних вершинах.

При этом в \cite{But24} рассматривалась более общая нестационарная управляемая система, заданная уравнениями
произвольного порядка нейтрального типа с негладкими коэффициентами, что, в свою очередь, потребовало
существенного развития понятия обобщенного решения краевой задачи для оптимальной траектории и в случае
интервала.

В частности, для ее постановки потребовалось введение семейства нелокальных квазипроизводных. Также было проведено их сравнение в локальном
случае $\tau=0$ с квазипроизводными, применяемыми для регуляризации сингулярных дифференциальных выражений с коэффициентами из пространств
обобщенных функций \cite{SavShk-99, MirShk-16, Bond-22}.

В \cite{But23-arXiv} показано, что в энергетический функционал на временн\'ом дереве целесообразно добавлять специальные веса, равные
вероятностям соответствующих сценариев.

В настоящей заметке обсуждается стохастическая интерпретация управляемой системы на дереве, основанная на
учете вероятностей всех возможных сценариев. А именно, к системе на конечном дереве приведет, например,
замена коэффициентов в уравнении~(\ref{1}) случайными процессами с дискретным временем и конечным числом
состояний.

Счетному числу состояний будет соответствовать более сложный и в то же время более общий случай бесконечного дерева, который и
рассматривается в данной работе.

\medskip
Для краткости ограничимся здесь управляемой системой, заданной уравнением
\begin{equation}\label{2}
y'(t)+by(t)=u(t), \quad 0<t<T, \quad T\in{\mathbb N},
\end{equation}
где $b\in{\mathbb C},$ которую выбором $u(t)\in L_2(0,T)$ требуется перевести из состояния
\begin{equation}\label{2-1}
y(0)=\varphi_0\in{\mathbb C}
\end{equation}
в состояние
\begin{equation}\label{2-2}
y(T)=\varphi_1\in{\mathbb C}.
\end{equation}
Для этого, очевидно, подойдет любое управление $u(t)=y'(t)+by(t),$ где $y(t)\in W_2^1[0,T]$ удовлетворяет краевым условиям (\ref{2-1}) и
(\ref{2-2}). Однако оптимальное управление, т.е. минимизирующее норму $\|u\|_{L_2(0,T)},$ будет единственно. При этом соответствующую
оптимальную траекторию $y(t)$ можно найти, решив краевую задачу
\begin{equation}\label{3}
-y''(t)-2i({\rm Im\,}b)y'(t)+|b|^2y(t)=0, \quad 0<t<T, \quad y(0)=\varphi_0, \quad y(T)=\varphi_1.
\end{equation}

Теперь предположим, что в моменты времени $t\in{\mathbb N}$ значение коэффициента $b$ может изменяться. А
именно, каждый раз он может принимать одно из счетного множества значений $\{\theta_l\}_{l\in{\mathbb
N}}\subset{\mathbb C},$ причем вероятность выпадения $b=\theta_l$ равна $p_l>0,$ и
$$
\sum_{l=1}^\infty p_l=1, \quad \sup_{l\in{\mathbb N}}{|\theta_l|}<\infty.
$$
В этом случае уравнение (\ref{2}) можно представить в виде системы дифференциальных уравнений с постоянными коэффициентами на некотором
бесконечном дереве ${\cal T}.$

В самом деле, при $t\in[0,1]$ функция $y_1(t):=y(t)$ является решением уравнения
\begin{equation}\label{4}
y_1'(t)+b_1y_1(t)=u_1(t), \quad 0<t<1,
\end{equation}
где через $b_1$ и $u_1(t)$ обозначены, соответственно, $b$ и $u(t)$ на интервале $(0,1),$ который, в свою очередь, параметризует первое
ребро $e_1=[v_0,v_1]$ нашего дерева ${\cal T}.$ Вершина~$v_0$ является корнем ${\cal T},$ а~$v_1$ соответствует первой точке $t=1$
ожидаемого изменения $b.$

Тогда при $t\in(1,2)$ будем иметь уже счетное число возможных уравнений
\begin{equation}\label{5}
\tilde y_j'(t)+b_j\tilde y_j(t)=u(t), \quad 1<t<2, \quad j\ge2,
\end{equation}
где $b_j=\theta_{j-1},$ а $\tilde y_j(t)$ обозначает соответствующее решение $y(t)$ на отрезке $[1,2].$

Уравнения в (\ref{5}) заданы, соответственно, на ребрах $e_j=[v_1,v_j],$ $j\ge2,$ образующих второй ярус дерева ${\cal T}.$ При этом в
вершине $v_1$ выполняются условия непрерывности
\begin{equation}\label{6}
y_1(1)=\tilde y_j(1),\quad j\ge2.
\end{equation}

В рамках рассматриваемой задачи оптимального управления естественно допустить в~(\ref{5}) вместо единой $u(t)$ различные правые части
$\tilde u_j(t)\in L_2(1,2),$ которые могут выбираться в зависимости от сценария $b=\theta_{j-1},$ реализовавшегося в момент времени $t=1.$

Обозначая $y_j(t):=\tilde y_j(t+1)$ и $u_j(t):=\tilde u_j(t+1)$ при $j\ge2,$ перепишем (\ref{4}) и (\ref{5}) в виде
\begin{equation}\label{7}
\ell_j y(t):=y_j'(t)+b_j y_j(t)=u_j(t), \quad 0<t<1, \quad j\in{\mathbb N}.
\end{equation}
Тогда (\ref{6}) становятся условиями склейки
\begin{equation}\label{8}
y_1(1)=y_j(0),\quad j\ge2.
\end{equation}
Вместе с начальным условием (\ref{2-1}), которое принимает вид
\begin{equation}\label{9}
y_1(0)=\varphi_0,
\end{equation}
уравнения (\ref{7}) и условия (\ref{8}) образуют задачу Коши на бесконечном графе-звезде.

Если $T=2,$ то построение дерева ${\cal T}$ завершено. В случае $T>2$ нужно повторить сделанные выше шаги. В
результате придем к задаче Коши на некотором бесконечном дереве ${\cal T},$ состоящем из $T$ ярусов. Однако
для ее записи удобно будет использовать другую нумерацию вершин, для которой мы введем обозначения
$$
{\mathbb N}_1:=\{2j-1:j\in{\mathbb N}\}, \quad {\mathbb N}_2:={\mathbb N}\setminus{\mathbb N}_1, \quad {\mathbb N}_\nu^0:={\mathbb
N}_\nu\cup\{0\}, \quad \nu=1,2.
$$
Именно, будем считать, что вершины $\{v_j\}_{j\in{\mathbb N}_1}$ являются внутренними, а
$\{v_j\}_{j\in{\mathbb N}_2^0}$ -- граничными, и $v_0$ по-прежнему обозначает корень ${\cal T}.$ При этом
всякое ребро имеет вид $e_j=[v_{k_j},v_j],$ $j\in{\mathbb N},$ где $k_j$ является отображением натуральных
чисел на ${\mathbb N}_1^0,$ задающим структуру ${\cal T}.$ А именно, для каждого $j\in{\mathbb N}_1^0$
множество $\{e_\nu\}_{\nu\in V_j},$ где $V_j:=\{\nu:k_\nu=j\},$ состоит из всех ребер, выходящих из вершины
$v_j.$ Очевидно, $\#V_0=1,$ поскольку корень является граничной вершиной. Пусть $k_1=0.$

Положим $k_j^{\{0\}}:=j$ и $k_j^{\{\nu+1\}}:=k_{k_j^{\{\nu\}}}$ при $\nu=0,\ldots,\nu_j,$ где $\nu_j$ таково, что $k_j^{\{\nu_j\}}=1.$
Тогда для $j\in{\mathbb N}$ цепочка ребер ${\cal E}_j:=\{e_{k_j^{\{\nu\}}}\}_{\nu=0,\ldots,\nu_j}$ образует единственный простой путь между
вершиной $v_j$ и корнем.

При $j\in{\mathbb N}_1$ множество $V_j$ является счетным. Обозначим его элементы через $\nu_{j,l},$ т.е. $\{\nu_{j,l}\}_{l\in{\mathbb
N}}=V_j$ при $j\in{\mathbb N}_1.$ Положим $b_{\nu_{j,l}}:=\theta_l$ и $\tilde p_{\nu_{j,l}}:=p_l$ при $l\in{\mathbb N}$ и $j\in{\mathbb
N}_1,$ а также
\begin{equation}\label{9-1}
\alpha_1:=1, \quad \alpha_j:=\prod_{\nu=0}^{\nu_j-1}\tilde p_{k_j^{\{\nu\}}}, \quad j\ge2.
\end{equation}
Другими словами, при $j\ge2$ будем иметь $\alpha_j=\tilde p_j\alpha_{k_j}.$

Таким образом, упорядоченные наборы $B_j:=[b_{k_j^{\{\nu\}}}]_{\nu=\nu_j,\ldots,0}$ при $j\in{\mathbb N}_2$ представляют собой все
реализации процесса $b$ на $(0,T),$ причем вероятность $B_j$ равна $\alpha_j.$

По-прежнему считаем, что каждое ребро $e_j$ параметризовано переменной $t\in[0,1],$ причем значение $t=0$ всегда соответствует его началу
$v_{k_j},$ а $t=1$ -- концу $v_j.$

Легко видеть, что таким образом будет охвачен и случай графа-звезды, если положить ${\mathbb N}_1=\{1\}$ с
сохранением остальных сделанных обозначений. Тогда, например, $V_1={\mathbb N}_2.$

На дереве ${\cal T}$ рассмотрим задачу Коши для системы уравнений (\ref{7}) (при новом определении $b_j)$ с начальным условием (\ref{9}) и
условиями склейки во внутренних вершинах
\begin{equation}\label{10}
y_{k_j}(1)=y_j(0), \quad j\ge2,
\end{equation}
которую можно интерпретировать, как разложение задачи Коши на интервале (\ref{2}), (\ref{2-1}) по всем
реализациям процесса~$b,$ причем вероятность $j$-го уравнения в (\ref{7}) равна $\alpha_j.$

Точнее, будем понимать задачу (\ref{2}), (\ref{2-1}) как задачу (\ref{7}), (\ref{9}), (\ref{10}) при $u_j(t)=u(t+k-1)$ для $e_j,$
образующих $k$-й ярус дерева ${\cal T},$ т.е. для которых $\#{\cal E}_j=k,$ где $k=1,\ldots,T.$

При $s+1\in{\mathbb N}$ рассмотрим весовое гильбертово пространство $W_{2,\alpha}^s$ со скалярным произведением
$\langle\,\cdot\,,\,\cdot\,\rangle_s$ и соответствующей нормой $\|\,\cdot\,\|_s=\sqrt{\langle\,\cdot\,,\,\cdot\,\rangle_s}:$
$$
W_{2,\alpha}^s=\prod_{j=1}^\infty W_2^s[0,1], \quad \langle y,\,z\rangle_s=\sum_{j=1}^\infty \alpha_j\langle y_j,\,z_j\rangle_{W_2^s[0,1]},
\quad y=[y_j],\,z=[z_j]\in W_{2,\alpha}^s.
$$
Нетрудно убедиться в справедливости следующего утверждения.

\begin{lemma}\label{lemm1}
Для всякого $u=[u_j]\in W_{2,\alpha}^0=:L_{2,\alpha}$ задача Коши (\ref{7}), (\ref{9}), (\ref{10}) имеет единственное решение $y=[y_j]\in
W_{2,\alpha}^1.$

С другой стороны, для любого $y=[y_j]\in W_{2,\alpha}^1$ справедливо $\ell y:=[\ell y_j]\in L_{2,\alpha}.$
\end{lemma}

Нашей целью будет являться нахождение такого управления $u\in L_{2,\alpha},$ которое приведет систему (\ref{7}), (\ref{9}), (\ref{10}) в
состояние
\begin{equation}\label{11}
y_j(1)=\varphi_1, \quad j\in{\mathbb N}_2,
\end{equation}
т.е. сразу для всех реализаций случайного процесса $b.$

Поскольку такое $u$ не единственно, будем искать его из условия минимума требуемых усилий $\|u\|_0.$ Приходим к вариационной задаче о
минимуме функционала энергии
\begin{equation}\label{12}
J(y):=\sum_{j=1}^\infty \alpha_j\int_0^1 |\ell_j y (t)|^2\,dt \to\min
\end{equation}
при условиях (\ref{9}), (\ref{10}), (\ref{11}), которую будем обозначать через ${\cal V}.$

\begin{theorem} \label{theo1}
Кортеж $y\in W_{2,\alpha}^1$ является решением вариационной задачи ${\cal V}$ тогда и только тогда, когда он обладает дополнительной
гладкостью $y\in W_{2,\alpha}^2$ и является решением краевой задачи ${\cal B}$ на дереве ${\cal T},$ состоящей из уравнений
\begin{equation}\label{13}
{\cal L}_jy(t):=-y_j''(t)-2i({\rm Im\,}b_j)y_j'(t)+|b_j|^2y_j(t)=0, \quad 0<t<1, \quad j\in{\mathbb N},
\end{equation}
и условий (\ref{9}), (\ref{10}), (\ref{11}) вместе с условиями типа Кирхгофа во внутренних вершинах
\begin{equation}\label{14}
y_j'(1)+\beta_jy_j(1)=\sum_{\nu\in V_j}\tilde p_\nu y_\nu'(0), \quad \beta_j:=b_j-\sum_{\nu\in V_j}\tilde p_\nu b_\nu, \quad j\in{\mathbb
N}_1.
\end{equation}
\end{theorem}

\begin{theorem} \label{theo2}
Оператор, порожденный дифференциальными выражениями (\ref{13}), с областью определения, состоящей из кортежей
$y=[y_j]\in W_{2,\alpha}^2,$ удовлетворяющих условиям (\ref{10}) и (\ref{14}), а также условиям (\ref{9}) и
(\ref{11}) при $\varphi_0=\varphi_1=0,$ является самосопряженным положительно определенным оператором в
$L_{2,\alpha}.$
\end{theorem}

Вместе с теоремой~\ref{theo1} следующая теорема дает однозначную разрешимость задачи ${\cal V}.$

\begin{theorem} \label{theo3}
Краевая задача ${\cal B}$ имеет единственное решение $y\in W_{2,\alpha}^2.$ Кроме того, существует такое $C,$ что выполняется априорная
оценка
$$
\|y\|_1\le C(|\varphi_0|+|\varphi_1|).
$$
\end{theorem}

Итак, установлены существование и единственность траектории $y=[y_j],$ являющейся оптимальной с учетом сразу
всех возможных реализаций процесса $b.$ Подставляя ее компоненты в уравнения (\ref{7}), получим
соответствующий набор управлений $u=[u_j],$ который должен применяться следующим образом. На интервале
$(0,1)$ в качестве управляющего воздействия выбирается первая его компонента $u(t):=u_1(t).$ Далее, если в
момент времени $t=1$ выясняется, что $b=\theta_l,$ то при $t\in(1,2)$ следует выбирать
$u(t):=u_{\nu_{1,l}}(t-1).$

Вообще, если к моменту времени $t=k\in\{1,\ldots,T-1\}$ реализовался набор сценариев $[b_{k_j^{\{\nu\}}}]_{\nu=k-1,\ldots,0}$ для
некоторого $j\in{\mathbb N}_1,$ а при $t=k$ выпадает $b=\theta_l,$ то на обозримую перспективу, т.е. при $t\in(k,k+1),$ выбирается
управление $u(t):=u_{\nu_{j,l}}(t-k).$

\medskip
{\sc Пример 1.} Пусть $b_j$ не зависят от $j,$ т.е. все $\theta_j$ равны первоначальному значению
коэффициента $b$ на интервале $(0,1).$ Тогда все пути ${\cal E}_j,$ $j\in{\mathbb N}_2,$ становятся
искусственными копиями единственной реализации $b.$ В силу (\ref{9-1}) будем иметь $\beta_j=0$ в (\ref{14})
для всех $j\in{\mathbb N}_1,$ а оптимальная траектория $y=[y_j]$ будет независимой от вероятностей $\{p_k\}.$

В самом деле, пусть ${\cal T}$ -- звезда. Тогда если в (\ref{14}) заменить любые два значения $\tilde
p_{j_1}$ и $\tilde p_{j_2}$ при $j_1,j_2\in {\mathbb N}_2(=V_1)$ их средним, то в силу образовавшейся
симметрии в таким образом модифицированной задаче ${\cal B}$ соответствующие компоненты $y_{j_1}(t)$ и
$y_{j_2}(t)$ решения $y$ равны. В противном случае оно было бы не единственным. Но это решение, очевидно,
будет решением и исходной задачи ${\cal B}.$ В силу произвольности выбора $j_1$ и $j_2,$ равны все $y_j(t)$
при $j\in{\mathbb N}_2.$ Аналогично компоненты решения задачи ${\cal B},$ соответствующие одному и тому же
ярусу произвольного дерева ${\cal T},$ также равны между собой. Чтобы в этом убедиться, достаточно будет
рассмотреть произвольное бинарное поддерево дерева ${\cal T}$ и на каждом его ярусе усреднить соответствующие
значения~$\tilde p_j.$

Таким образом, краевая задача ${\cal B}$ оказывается эквивалентной задаче (\ref{3}) на интервале, решение
которой будет связано с решением ${\cal B}$ следующим образом:
$$
y(t+k-1)=y_{k_2^{\{T-k\}}}(t), \quad 0\le t\le1, \quad k=1,\ldots,T.
$$

\begin{remark}
Тогда как рассмотренный процесс $b$ представляет собой набор независимых в совокупности случайных величин,
конструкция функционала (\ref{12}) не зависит от этого ограничения. В самом деле, каждая величина $\tilde
p_s,$ $s\in V_j,$ фигурирующая в определении (\ref{9-1}), вообще говоря, может быть условной вероятностью
сценария, отвечающего ребру $e_s,$ при условии, что к моменту времени $t=\nu_j+1$ реализовался набор ребер
${\cal E}_j.$

Последнее наблюдение позволяет сделать наборы $\{\tilde p_{\nu_{j,l}}\}_{l\in{\mathbb N}}$ зависящими от
$j\in{\mathbb N}_1.$ При этом каким-либо $\tilde p_{\nu_{j,l}}$ можно позволять обращаться в нуль, что будет
означать отсутствие соответствующих ребер $e_{\nu_{j,l}}.$ Поэтому приведенная конструкция ${\cal T}$
распространяется на случай наличия вершин конечной кратности и, в частности, конечного дерева.

Таким образом, теоремы~\ref{theo1}--\ref{theo3} останутся справедливы в случае любого дискретного случайного
процесса $b$ с дискретным временем ${\mathbb N}$ и ограниченным множеством состояний в ${\mathbb C}.$ При
этом для конечного дерева теоремы \ref{theo1} и \ref{theo3} являются, соответственно, частными случаями
теорем~5 и~6 из \cite{But24}.

Кроме того, теоремы~\ref{theo1} и~\ref{theo3} будут также справедливы, если заменить условия (\ref{11}) на условия
$y_j(1)=\psi_j\in{\mathbb C},$ $j\in{\mathbb N}_2,$ требуя только $\varphi_1:=\sup|\psi_j|<\infty.$

Наконец, отметим, что приведенная конструкция дерева легко распространяется на случай, когда длины ребер могут различаться, оставаясь
ограниченными в совокупности. В частности, это позволяет охватить и неравномерно распределенное дискретное время.
\end{remark}

\end{fulltext}

\end{document}